# Numbers in the base $e^\pi$

Simon Plouffe
September 19, 2025

Summary

A large-scale experiment was conducted to find formulas relating to the base $e^\pi$. The numbers in this base are

$$x = \sum_{n=0}^{\infty} \frac{a(n)}{e^{\pi n}}$$

Where $a(n)$ is taken from the OEIS catalog. These experiments were inspired by several facts. Indeed, it is known that the formula generating the partitions of integers is generated by

$$\prod_{k \geq 1}^{\infty} \frac{1}{1-x^k} = \sum_{n=0}^{\infty} p(n)\, x^n$$

is equal to

$$\frac{2^{3/8}\, \Gamma(3/4)}{\pi^{1/4}\, e^{\pi/24}} \tag{1}$$

when evaluated at the point $x = e^{-\pi}$. By analyzing the 387500 sequences of the OEIS catalog, the model that was used is based on the fact that the infinite sum evaluated at $e^\pi$, is an expression that can be detected using a program like *lindep* from Pari-Gp. The process made it possible to find 793 expressions similar to (1).

Most of the known real numbers in base 10 that come from classical mathematical analysis have no pattern in their decimal expansion. On the other hand, in base $e^{\pi}$ it is much richer since a class of numbers have coefficients that come precisely from analysis, combinatorics and number theory. The examples are numerous and varied.

To detect a relationship with integers, it is necessary to take the logarithm of the sum and evaluate to a certain precision (160 decimal digits) whether the '=' sign is verified.

$$\log\left(\sum_{n=0}^{\infty} \frac{a(n)}{e^{k\pi n}}\right) \equiv [\pi, \log(\pi), \log(2), \log(3), \log \Gamma(3/4)] \qquad (2)$$

The sign $\equiv$ means that there is an identity with 0 or in other words, the logarithm of the sum is a linear combination of the list of constants. This assumption seems to be correct since more than 793 expressions have been identified. The first version of this article published in 2023 had identified 659. By adding elements to the list (2):

$$\pi, \log(\pi), \log(2), \log(3), \log(5), \log \Gamma(1/3), \log \Gamma(1/4),$$
$$\log \Gamma(1/5), \log \Gamma(2/5), \log \Gamma(1/6), \log \Gamma(1/8), \log \Gamma(3/8), \log \Gamma(1/10),$$
$$\log \Gamma(3/10), \log \Gamma(1/12), \log \Gamma(5/12).$$

We obtain more identities, from 659 to 793 including those related to the most well-known Ramanujan functions.

The remarkable thing about these expressions is the fact that despite the great variety of combinatorial contexts, they all have the same pattern if we evaluate the sequence at the point $e^{\pi}$ The experiment is conclusive, since 363 sequences refer to theta functions, 262 to Ramanujan functions and variants and 57 to partitions (sharings) of integers. However, work remains to be done since there are thousands of other sequences referring to theta functions and thousands to partitions in all contexts.

Note: The extended sequences come from the 'b' files on the OEIS site, normally a sequence has about 3 full lines of about 80 characters of terms. In many cases it has long been thought that an extension of the basic sequence was necessary. The choice was made to take the first 2000 terms of each sequence (when available) and with a precision of 160 decimals for the sums evaluated.

Index: Most sequences are defined from 0, some from 1 or another starting point. This doesn't actually change the final result since the factor $e^{\pi}$ then appears in the final expression, it's just an offset. So, I took this point as a starting point. In some cases the exponent is even fractional.

Each page contains the sequence number, name, formula, first few terms of the sequence and value in $x = e^{-\pi}$.

## Other approaches

### Basic expansion in the base $e^\pi$

We can go the other way from a known value. A good example is the sequence A000122 which lists the number of ways to represent n as a sum of 1 square. Evaluated at $x = e^{-\pi}$ it gives us the value $\frac{\pi^{1/4}}{\Gamma\left(\frac{3}{4}\right)}$ but if we perform the expansion in base $e^\pi$ we get the sequence:

$$1,2,0,0,2,0,0,0,0,2,0,0,0,0,0,0,2,0,0,0,0,0,0,0,0,2,0, \dots$$

It works of course if we know the number in advance $\frac{\pi^{1/4}}{\Gamma\left(\frac{3}{4}\right)}$.

This method only works if the terms of the sequence $a(n)$ do not exceed the value $e^\pi$. The expansion of a number in base $e^\pi$ uses the same process that can be used to expand a real number in base 2, 10 or 16. It is therefore possible from values similar to formula (1) to go backwards. A program has found some examples of this but it is not enough to make it a method.

### The Rogers-Ramanujan identities

Here we consider the two functions $G(q)$ and $H(q)$.

$$G(q) = \sum_{n=0}^{\infty} \frac{q^{n^2}}{(q;q)_n} = \frac{1}{(q;q^5)_\infty (q^4;q^5)_\infty} = 1 + q + q^2 + q^3 + 2q^4 + 2q^5 + 3q^6$$

$$H(q) = \sum_{n=0}^{\infty} \frac{q^{n^2+n}}{(q;q)_n} = \frac{1}{(q^2;q^5)_\infty (q^3;q^5)_\infty} = 1 + q^2 + q^3 + q^4 + q^5 + 2q^6 + \cdots$$

The 2 series are the 2 sequences: A003114 and A003105 of the OEIS.

And we are interested in the $R(q)$ the well-known function

$$R(q) = q^{1/5} \prod_{k=0}^{\infty} \frac{(1-q^{5k+1})(1-q^{5k+4})}{(1-q^{5k+2})(1-q^{5k+3})} = q^{1/5} \frac{H(q)}{G(q)}$$

Which has a remarkable continued fraction expansion.

$$R(q) = q^{1/5} \left[ 1 + \cfrac{q}{1 + \cfrac{q^2}{1 + \cfrac{q^3}{1 + \cdots}}} \right]$$

When $q = e^{-\pi}$ the value is algebraic and the numeric value is 0.5114284554037...,

Or more precisely

$$\frac{1}{4}(\sqrt{5}+1)(\sqrt{5}-\sqrt{\sqrt{5}+2})(\sqrt{\sqrt{5}+2}+\sqrt[4]{5})$$

But what is less so is the value when it is evaluated at $e^{-\pi/n}$

| n | value at $e^{-\pi/n}$ |
|---|---|
| 1 | .51142845540370351929463301354268 |
| 2 | .60900295189734706892868001337065 |
| 3 | .61729972602946859110062866515646 |
| 4 | .61797449361217335219105975903425 |
| 5 | .61802916937167616162641498260454 |
| 6 | .61803359836699331081163334071515 |
| 7 | .61803395712786281601119176526411 |
| 8 | .61803398618842825100587467815680 |
| 9 | .61803398854240941388707852692804 |
| 10 | .61803398873087990188418822923 94 |
| 11 | .61803398748533449139237316036 80 |
| 12 | .61803398749784571351332060666 00 |
| 13 | .61803398749885915493618176183 33 |
| 14 | .61803398749894124631843375359 11 |
| 15 | .61803398749894789593318999351 78 |
| 16 | .61803398749894843456922210066 10 |
| 17 | .61803398749894847820013693101 02 |
| 18 | .61803398749894848173435413048 87 |
| 19 | .61803398749894848202063488489 60 |
| 20 | .61803398749894848204382436808 79 |
| 21 | .61803398749894848204570277633 69 |
| 22 | .61803398749894848204585493227 31 |
| 23 | .61803398749894848204586725730 31 |
| 24 | .61803398749894848204586825566 26 |

The values at m are almost exactly the n'th decimal of $\phi - 1$. The known limit of $\phi - 1$ is easy to check if we examine the continued fraction of $R(q)$.

The same approximation phenomenon occurs with A000122, $a(n)$ is the number of solutions of $k^2 = n$. The sequence reads: 1, 2, 0, 0, 2, 0, 0, 0, 0, 2, 0, 0, 0, 0, 0, 0, 2, 0, 0, …
If we evaluate in $e^{-\pi/n}$ we see another pattern appear by constructing a small table of values we have

| n | value at $e^{-\pi/n}$ |
|---|---|
| 1 | 1.08643481121330801457531612151 02 |
| 2 | 1.41949548808376612336218673135 17 |
| 3 | 1.73233035889803357021221785709 57 |
| 4 | 2.00001394936942483598255871491 14 |
| 5 | 2.23606865145840390415024699212 08 |
| 6 | 2.44948977468735155440609102979 78 |

```
 7 2.6457513110645905905016157536392
 8 2.8284271247461900976033774484194
 9 3.0000000000031532911056038691332
10 3.1622776601852296904234046087410
11 3.3166247903554063591967389627611
12 3.4641016151377548808904659455994
13 3.6055512754639893063354964239870
14 3.7416573867394138617643526584660
15 3.8729833462074168852057766566470
16 4.0000000000000000000011832276929
17 4.1231056256176605498214625615352
```

We can quickly guess that the function $\sqrt{n}$ to a precision of n decimal places. The precision increases when n >> 1.

## The calculation of $\Gamma(3/4)$ or can this number be calculated with great precision?

The starting point is one of the values of the Jacobi function $\vartheta_3(q)$ which when evaluated at $q = e^{-\pi}$, the infinite sum is then

$$\vartheta_3(q) = \sum_{n=-\infty}^{\infty} q^{n^2} = 1 + 2q + 2q^4 + 2q^9 + 2q^{16} + \cdots$$

Is equal to $\frac{\pi^{1/4}}{\Gamma(3/4)}$, this series is in fact the sequence A000122 of the OEIS catalog.

As for the $\vartheta_3(q)$ Jacobi function, it has an interesting property. For values of $m \in \mathbb{N}^*$, we have

$$\vartheta_3(e^{-m\pi}) = \frac{A_m \pi^{\frac{1}{4}}}{\Gamma\left(\frac{3}{4}\right)}$$

With $A_m$ a possibly high degree algebraic number. The values of $A_m$ are known for some m that we have managed [9] to calculate explicitly.
As here with m = 17 we have [8].

$$\vartheta_3(e^{-17\pi}) = \frac{\sqrt[4]{\pi}}{\Gamma\left(\frac{3}{4}\right)} \frac{\sqrt{2}\left(1 + \sqrt[4]{17}\right) + \sqrt[8]{17}\sqrt{5 + \sqrt{17}}}{\sqrt{17 + 17\sqrt{17}}}$$

The complexity of the equation poses a problem of organizing the calculation into separate steps. So, the idea is to find the highest possible m and from there use the fast convergence properties of the function. But there remains the problem of the exact representation of $\vartheta_3(e^{-m\pi})$.

We manage to find a high value with m = 289. Using Pari-GP is the *algdep* with a precision of 32000 decimal digits, an algebraic number is detected and is of degree 136. It would be quite difficult, if not impossible, to try to find the expression in radicals even if it exists. Another problem then arises. If an equation in radicals is found, it is still necessary to be able to evaluate it with great precision.

The solution that is by far the simplest is to use Newton's algorithm. Starting from the first approximation, for example, at 32,000 decimal digits, it only takes about thirty iterations to have the equivalent of thousands of billions of decimal places, the precision doubling at each step.

In summary: to evaluate $\Gamma\left(\frac{3}{4}\right)$ very precisely we separate the calculation into 3 parts.
1- Calculating $e^\pi$, a precise value $\pi$ is needed at the start, this is the most difficult part.
2- Calculation of $\sqrt{\sqrt{\pi}}$ by Newton's method.
3- We assume that $A_m$ exists and we test with different values of m that satisfy an algebraic equation, the values of m found are m = 2, 3, 4, 5, 6, 8, 9, 13, 16, 17, 21, 25, 29, 33, 36, 37, 41, 49, 53, 61, 64, 65, 73, 85, 89, …. We suspect that when m is a perfect square the polynomials are more easily detectable. We then push the 'algdep' program to the maximum with 32000 decimal digits and an algebraic equation that is as small as possible.
4- One was found: 289, so with 1806 terms of the series $\vartheta_3(e^{-289\pi})$ we would have a precision of 1 billion digits and 2305 billion digits if we take only 100000 terms.
5- Once m = 289 is found, it is enough to apply Newton with his formula where $x_0$ is the root of the polynomial.

$$x_{k+1} = x_k - \frac{f(x_k)}{f'(x_k)}$$

To find the real roots, in this case there are only 2. $f(x)$ is the polynomial and its derivative $f'(x)$. Once steps 1-2 are done, we isolate $\Gamma\left(\frac{3}{4}\right)$ and the calculation is complete.

You can consult the list of polynomials found here:
https://plouffe.fr/articles/polynomials%20for%20%20theta_3.pdf

According to Alexander Yee, see [10], the calculation of $e^\pi$ is by far the slowest, the other 2 calculations with Newton's method do not pose any problem, the algorithm has already been proven. The calculation of $f(x)$ with an m = 289 poses a bulge problem. The bulge comes from the fact that the coefficients of the polynomial are enormous, the coefficient of $x^{136}$ is simple but it has 156 digits, it is $17^{136}$, which says bulge says slowness. We could gain speed by taking a much smaller value of m like 5 and have an algebraic number very easy to evaluate with great precision but there remains the bottleneck of the value of $e^{\pi}$ to be carried out.

## Appendix: first pages and last pages of the document.

**A000009** Expansion of Product_{m >= 1} (1 + x^m); number of partitions of n into distinct parts; number of partitions of n into odd parts.

$$\frac{e^{\frac{\pi}{24}} 2^{7/8}}{2}$$

1.045250214354711942547595012203­6

1,1,1,2,2,3,4,5,6,8,10,12,15,18,22,27,32,38,46,54,64,76,89,104,122,142,165,192,222,
256,296,340,390,448,512,585,668,760,864,982,1113,1260,1426,1610,1816,2048,2304,
2590,2910,3264,3658,4097,4582,5120,5718,6378

---

**A000041** a(n) is the number of partitions of n (the partition numbers).

$$\frac{e^{-\frac{\pi}{24}} 2^{3/8} \Gamma\left(\frac{3}{4}\right)}{\pi^{1/4}}$$

1.0472094700460421302197988631655

1,1,2,3,5,7,11,15,22,30,42,56,77,101,135,176,231,297,385,490,627,792,1002,1255,1575,
1958,2436,3010,3718,4565,5604,6842,8349,10143,12310,14883,17977,21637,26015,
31185,37338,44583,53174,63261,75175,89134,105558,124754,147273,173525

---

**A000118** Number of ways of writing n as a sum of 4 squares; also theta series of four-dimensional cubic lattice Z^4.

$$\frac{\pi}{\Gamma\left(\frac{3}{4}\right)^4}$$

1.3932039296856768591842462603253

1,8,24,32,24,48,96,64,24,104,144,96,96,112,192,192,24,144,312,160,144,256,288,192,
96,248,336,320,192,240,576,256,24,384,432,384,312,304,480,448,144,336,768,352,
288,624,576,384,96,456,744,576,336,432,960,576,192

---

**A000122** Expansion of Jacobi theta function theta_3(x) = Sum_{m =-oo..oo} x^(m^2) (number of integer solutions to k^2 = n).

$$\frac{\pi^{1/4}}{\Gamma\left(\frac{3}{4}\right)}$$

1.0864348112133080145753161215102

1,2,0,0,2,0,0,0,0,2,0,0,0,0,0,0,2,0,0,0,0,0,0,0,2,0,0,0,0,0,0,0,0,0,2,0,0,0,0,0,0,0,0,

0,0,0,2,0,0,0,0,0,0,0,0,0,0,0,0,0,2,0,0,0,0,0,0,0,0,0,0,0,0,0,0,0,2,0,0,0,0,0,0,0,0,0,0,0,0,0,0,0,0,0,0,2,0,0,0,0

---

### A000132 Number of ways of writing n as a sum of 5 squares.

$$\frac{\pi^{5/4}}{\Gamma\left(\frac{3}{4}\right)^5}$$

1.51362524832969719205469025639191

1,10,40,80,90,112,240,320,200,250,560,560,400,560,800,960,730,480,1240,1520,752, 1120,1840,1600,1200,1210,2000,2240,1600,1680,2720,3200,1480,1440,3680,3040, 2250,2800,3280,4160,2800,1920,4320,5040,2800,3472,5920

---

### A000141 Number of ways of writing n as a sum of 6 squares.

$$\frac{\pi^{3/2}}{\Gamma\left(\frac{3}{4}\right)^6}$$

1.64445516091677103106944564291411

1,12,60,160,252,312,544,960,1020,876,1560,2400,2080,2040,3264,4160,4092,3480, 4380,7200,6552,4608,8160,10560,8224,7812,10200,13120,12480,10104,14144,19200, 16380,11520,17400,24960,18396,16440,24480,27200

---

### A000143 Number of ways of writing n as a sum of 8 squares.

$$\frac{\pi^{2}}{\Gamma\left(\frac{3}{4}\right)^8}$$

1.94101718969161242994989600479831

1,16,112,448,1136,2016,3136,5504,9328,12112,14112,21312,31808,35168,38528,56448, 74864,78624,84784,109760,143136,154112,149184,194688,261184,252016,246176, 327040,390784,390240,395136,476672,599152,596736,550368,693504,859952

---

### A000144 Number of ways of writing n as a sum of 10 squares.

$$\frac{\pi^{5/2}}{\Gamma\left(\frac{3}{4}\right)^{10}}$$

2.29106139238113749228606857623721

1,20,180,960,3380,8424,16320,28800,52020,88660,129064,175680,262080,386920, 489600,600960,840500,1137960,1330420,1563840,2050344,2611200,2986560,

3358080,4194240,5318268,5878440,6299520,7862400,9619560

---

*A000145 Number of ways of writing n as a sum of 12 squares.*

$$\frac{\pi^3}{\Gamma\left(\frac{3}{4}\right)^{12}}$$

2.70423277626580330600018095350677

1,24,264,1760,7944,25872,64416,133056,253704,472760,825264,1297056,1938336,
2963664,4437312,6091584,8118024,11368368,15653352,19822176,24832944,
32826112,42517728,51425088,61903776,78146664,98021616

---

*A000152 Number of ways of writing n as a sum of 16 squares.*

$$\frac{\pi^4}{\Gamma\left(\frac{3}{4}\right)^{16}}$$

3.76754773067832495079594086144133

1,32,480,4480,29152,140736,525952,1580800,3994080,8945824,18626112,36714624,
67978880,118156480,197120256,321692928,509145568,772845120,1143441760,
1681379200,2428524096,3392205824,4658843520,6411152640

---

*A000156 Number of ways of writing n as a sum of 24 squares.*

$$\frac{\pi^6}{\Gamma\left(\frac{3}{4}\right)^{24}}$$

7.31287490823025420018915311755607

1,48,1104,16192,170064,1362336,8662720,44981376,195082320,721175536,
2319457632,6631997376,17231109824,41469483552,93703589760,200343312768,
407488018512,793229226336,1487286966928,2697825744960,4744779429216

---

*A000594 Ramanujan's tau function (or Ramanujan numbers, or tau numbers).*

$$\frac{e^\pi \pi^6}{512\,\Gamma\left(\frac{3}{4}\right)^{24}}$$

0.33051755959632854743859758293250

1,-24,252,-1472,4830,-6048,-16744,84480,-113643,-115920,534612,-370944,-577738,
401856,1217160,987136,-6905934,2727432,10661420,-7109760,-4219488,-12830688,

18643272,21288960,-25499225,13865712,-73279080,24647168

---

A000700 Expansion of Product_{k>=0} (1 + x^(2k+1)); number of partitions of n into distinct odd parts; number of self-conjugate partitions; number of symmetric Ferrers graphs with n nodes.

$$\frac{e^{-\frac{\pi}{24}} 2^{1/4}}{}$$

1.0432982626446870125278756888156

1,1,0,1,1,1,1,1,2,2,2,2,3,3,3,4,5,5,5,6,7,8,8,9,11,12,12,14,16,17,18,20,23,25,26,29,33,35, 37,41,46,49,52,57,63,68,72,78,87,93,98,107,117,125,133,144,157,168,178,192,209, 223,236,255,276,294,312,335,361,385

---

A000712 Generating function = Product_{m>=1} 1/(1 - x^m)^2; a(n) = number of partitions of n into parts of 2 kinds.

$$\frac{e^{-\frac{\pi}{12}} 2^{3/4} \Gamma\left(\frac{3}{4}\right)^2}{\sqrt{\pi}}$$

1.0966476741541124095724129797593

1,2,5,10,20,36,65,110,185,300,481,752,1165,1770,2665,3956,5822,8470,12230,17490, 24842,35002,49010,68150,94235,129512,177087,240840,326015,439190,589128, 786814,1046705,1386930,1831065,2408658,3157789,4126070,5374390

---

A000716 Number of partitions of n into parts of 3 kinds.

$$\frac{2 e^{-\frac{\pi}{8}} 2^{1/8} \Gamma\left(\frac{3}{4}\right)^3}{\pi^{3/4}}$$

1.1484198296781527497673927188231

1,3,9,22,51,108,221,429,810,1479,2640,4599,7868,13209,21843,35581,57222,90882, 142769,221910,341649,521196,788460,1183221,1762462,2606604,3829437,5590110, 8111346,11701998,16790136,23964594,34034391,48104069,67679109,94800537, 132230021,183686994,254170332

---

A000727 Expansion of Product_{k >= 1} (1 - x^k)^4.

$$\frac{e^{\frac{\pi}{6}} \pi \sqrt{2}}{4 \Gamma\left(\frac{3}{4}\right)^4}$$

0.83150670626724744966596392634565

1,-4,2,8,-5,-4,-10,8,9,0,14,-16,-10,-4,0,-8,14,20,2,0,-11,20,-32,-16,0,-4,14,8,-9,20,26,0,
2,-28,0,-16,16,-28,-22,0,14,16,0,40,0,-28,26,32,-17,0,-32,-16,-22,0,-10,32,-34,-8,14,
0,45,-4,38,8,0,0,-34,-8,38,0,-22,-56,2,-28,0,0,-10,20,64,-40,-20,44

---

### A000728 Expansion of Product_{n>=1} (1-x^n)^5.

$$\frac{e^{\frac{5\pi}{24}} \pi^{5/4} 2^{1/8}}{4 \Gamma\left(\frac{3}{4}\right)^5}$$

0.79402137781535628587400492096605

1,-5,5,10,-15,-6,-5,25,15,-20,9,-45,-5,25,20,10,15,20,-50,-35,-30,55,-50,15,80,1,50,-35,
-45,-15,5,-50,-25,-55,85,51,50,10,-40,65,10,-10,-115,50,-115,-100,85,80,-30,5,20,45,
70,65,45,-55,-100

---

### A000729 Expansion of Product_{k >= 1} (1 - x^k)^6.

$$\frac{e^{\frac{\pi}{4}} \pi^{3/2} 2^{3/4}}{8 \Gamma\left(\frac{3}{4}\right)^6}$$

0.75822593332778584280344454334574

1,-6,9,10,-30,0,11,42,0,-70,18,-54,49,90,0,-22,-60,0,-110,0,81,180,-78,0,130,-198,0,
-182,-30,90,121,84,0,0,210,0,-252,-102,-270,170,0,0,-69,330,0,-38,420,0,-190,-390,0,
-108,0,0,0,-300,99,442,210,0,418,-294,0,0,-510,378,-540,138,0

---

### A000730 Expansion of Product_{n>=1} (1 - x^n)^7.

$$\frac{e^{\frac{7\pi}{24}} \pi^{7/4} 2^{3/8}}{8 \Gamma\left(\frac{3}{4}\right)^7}$$

0.72404419079064411620167339982762

1,-7,14,7,-49,21,35,41,-49,-133,98,-21,126,112,-176,-105,-126,140,-35,147,259,98,-420,
-224,238,-455,273,-14,322,406,-35,-7,-637,-196,245,-181,-574,462,147,924,217,-329,
-140,-7,-371,-777

---

### A000731 Expansion of Product (1 - x^k)^8 in powers of x.

$$\frac{e^{\frac{\pi}{3}} \pi^2}{8\, \Gamma\left(\frac{3}{4}\right)^8}$$

0.691403402567406529188714128773392

*1,-8,20,0,-70,64,56,0,-125,-160,308,0,110,0,-520,0,57,560,0,0,182,-512,-880,0,1190,
-448,884,0,0,0,-1400,0,-1330,1000,1820,0,-646,1280,0,0,-1331,-2464,380,0,1120,0,
2576,0,0,-880,1748,0,-3850,0,-3400,0,2703,4160,-2500,0,3458*

---

*A000735 Expansion of Product_{k>=1} (1 - x^k)^12.*

$$\frac{e^{\frac{\pi}{2}} \pi^3 \sqrt{2}}{32\, \Gamma\left(\frac{3}{4}\right)^{12}}$$

0.574906565970791942075109707726666

*1,-12,54,-88,-99,540,-418,-648,594,836,1056,-4104,-209,4104,-594,4256,-6480,-4752,
-298,5016,17226,-12100,-5346,-1296,-9063,-7128,19494,29160,-10032,-7668,-34738,
8712,-22572,21812,49248,-46872,67562,2508,-47520,-76912,-25191,67716*

---

*A000739 Expansion of Product_{k>=1} (1 - x^k)^16.*

$$\frac{e^{\frac{2\pi}{3}} \pi^4}{64\, \Gamma\left(\frac{3}{4}\right)^{16}}$$

0.478038665081787213518128666645550

*1,-16,104,-320,260,1248,-3712,1664,6890,-7280,-5568,-4160,33176,4640,-74240,29824,
14035,54288,27040,-142720,1508,-110240,289536,222720,-380770,-83200,-123904,
142912,7640,408000,386048*

---

*A001934 Expansion of 1/theta_4(q)^2 in powers of q.*

$$\frac{\sqrt{2}\, \Gamma\left(\frac{3}{4}\right)^2}{\sqrt{\pi}}$$

1.19814023473559220743992249228804

*1,4,12,32,76,168,352,704,1356,2532,4600,8160,14176,24168,40512,66880,108876,
174984,277932,436640,679032,1046016,1597088,2418240,3632992,5417708,
8022840,11802176,17252928,25070568,36223424,52053760,74414412*

---

A001935 Number of partitions with no even part repeated; partitions of n in which no parts are multiples of 4.

$$\frac{e^{\frac{\pi}{8}}\sqrt{2}}{2}$$

1.0472058180553657181251543939073

1,1,2,3,4,6,9,12,16,22,29,38,50,64,82,105,132,166,208,258,320,395,484,592,722,876, 1060,1280,1539,1846,2210,2636,3138,3728,4416,5222,6163,7256,8528,10006,11716, 13696,15986,18624,21666,25169,29190,33808,39104,45164

---

A001936 Expansion of q^(-1/4) * (eta(q^4) / eta(q))^2 in powers of q.

$$\frac{e^{\frac{\pi}{4}}}{2}$$

1.0966400253690077282798848296394

1,2,5,10,18,32,55,90,144,226,346,522,777,1138,1648,2362,3348,4704,6554,9056,12425, 16932,22922,30848,41282,54946,72768,95914,125842,164402,213901,277204, 357904,460448,590330,754368,960948,1220370,1545306

---

A001937 Expansion of (psi(x^2) / psi(-x))^3 in powers of x where psi() is a Ramanujan theta function.

$$\frac{e^{\frac{3\pi}{8}}\sqrt{2}}{4}$$

1.1484078148788087524710011715016

1,3,9,22,48,99,194,363,657,1155,1977,3312,5443,8787,13968,21894,33873,51795, 78345,117312,174033,255945,373353,540486,776848,1109040,1573209,2218198, 3109713,4335840,6014123,8300811,11402928,15593702,21232521,28790667, 38884082

---

A001938 Expansion of k/(4*q^(1/2)) in powers of q, where k defined by sqrt(k) = theta_2(0, q)/theta_3(0, q).

$$\frac{e^{\frac{\pi}{2}}\sqrt{2}}{8}$$

0.85038029420627578205997577522294

1,-4,14,-40,101,-236,518,-1080,2162,-4180,7840,-14328,25591,-44776,76918,-129952, 216240,-354864,574958,-920600,1457946,-2285452,3548550,-5460592,8332425,

-12614088,18953310,-28276968,41904208,-61702876,90304598,-131399624

---

A001939 Expansion of (psi(-x) / phi(-x))^5 in powers of x where phi(), psi() are Ramanujan theta functions.

$$\frac{e^{\frac{5\pi}{8}} \sqrt{2}}{8}$$

1.2593899752426635611873132752101

1,5,20,65,185,481,1165,2665,5820,12220,24802,48880,93865,176125,323685,583798,
1035060,1806600,3108085,5276305,8846884,14663645,24044285,39029560,
62755345,100004806,158022900,247710570,385366265,595212280,913040649,
1391449780

---

A001940 Absolute value of coefficients of an elliptic function.

$$\frac{e^{\frac{3\pi}{4}}}{8}$$

1.3188405092747202734708839309306

1,6,27,98,309,882,2330,5784,13644,30826,67107,141444,289746,578646,1129527,
2159774,4052721,7474806,13569463,24274716,42838245,74644794,128533884,
218881098,368859591,615513678,1017596115,1667593666,2710062756,4369417452

---

A001941 Absolute values of coefficients of an elliptic function.

$$\frac{e^{\frac{7\pi}{8}} \sqrt{2}}{16}$$

1.3810974543995885825897999378618

1,7,35,140,483,1498,4277,11425,28889,69734,161735,362271,786877,1662927,3428770,
6913760,13660346,26492361,50504755,94766875,175221109,319564227,575387295,
1023624280,1800577849,3133695747,5399228149,9214458260,15584195428

---

A002107 Expansion of Product_{k>=1} (1 - x^k)^2.

$$\frac{e^{\frac{\pi}{12}} \sqrt{\pi}\, 2^{1/4}}{2\, \Gamma\!\left(\frac{3}{4}\right)^2}$$

0.91186989547152364013707755048620

1,-2,-1,2,1,2,-2,0,-2,-2,1,0,0,2,3,-2,2,0,0,-2,-2,0,0,-2,-1,0,2,2,-2,2,1,2,0,2,-2,-2,2,0,-2,0,

*A319307 Expansion of theta_4(q)^16 in powers of q = exp(Pi i t).*

$$\frac{\pi^4}{16\,\Gamma\!\left(\frac{3}{4}\right)^{16}}$$

0.23547173316795309424746303838833

1,-32,480,-4480,29152,-140736,525952,-1580800,3994080,-8945824,18626112,
-36714624,67978880,-118156480,197120256,-321692928,509145568,-772845120,
1143441760,-1681379200,2428524096,-3392205824,4658843520,-6411152640,
8705492608,-11488092896

---

*A319308 Expansion of theta_4(q)^20 in powers of q = exp(Pi i t).*

$$\frac{\pi^5}{32\,\Gamma\!\left(\frac{3}{4}\right)^{20}}$$

0.16403007198935613910580890178095

1,-40,760,-9120,77560,-497648,2508000,-10232640,34729720,-100906760,259114704,
-606957280,1327461600,-2738111280,5341699520,-9915552192,17701924600,
-30615844560,51294999960,-83279292960,131880275664,-204949382400,
312126610080,-464844224960,680432137440

---

*A319309 Expansion of theta_4(q)^24 in powers of q = exp(Pi i t).*

$$\frac{\pi^6}{64\,\Gamma\!\left(\frac{3}{4}\right)^{24}}$$

0.11426367044109772187795551746189

1,-48,1104,-16192,170064,-1362336,8662720,-44981376,195082320,-721175536,
2319457632,-6631997376,17231109824,-41469483552,93703589760,-200343312768,
407488018512,-793229226336,1487286966928,-2697825744960,4744779429216,
-8110465650176

---

*A319310 Expansion of theta_4(q)^28 in powers of q = exp(Pi i t).*

$$\frac{\pi^7}{128\,\Gamma\!\left(\frac{3}{4}\right)^{28}}$$

0.07959629733942323193011628491713

1,-56,1512,-26208,327656,-3147984,24189984,-152867520,811401192,-3681079640,
14500933104,-50376047904,156797510688,-444306558864,1163495873088,

-2851049839680,6597606440936,-14512424533488,30505974273096,
-61591664700384,119983597365744,-226303038736128

---

A319552 Expansion of 1/theta_4(q)^3 in powers of q = exp(Pi i t).

$$\frac{2^{3/4}\,\Gamma\!\left(\dfrac{3}{4}\right)^{3}}{\pi^{3/4}}$$

1.31147941617165978854278160756938

1,6,24,80,234,624,1552,3648,8184,17654,36816,74544,147056,283440,535008,990912,
1803882,3232224,5707624,9943536,17106960,29088352,48922320,81438528,
134261584,219336630,355242288,570675904,909674688,1439394192,2261635168,
3529838208

---

A319553 Expansion of 1/theta_4(q)^8 in powers of q = exp(Pi i t).

$$\frac{4\,\Gamma\!\left(\dfrac{3}{4}\right)^{8}}{\pi^{2}}$$

2.06077515502864631635689459626580

1,16,144,960,5264,25056,106944,418176,1520784,5201232,16871648,52252992,
155341248,445226848,1234726272,3323392128,8704504976,22234655520,
55498917840,135595345600,324759439584,763505859072,1764050361152,
4009763323008,8975341703616,19800832628336

---

A319554 Expansion of 1/theta_4(q)^12 in powers of q = exp(Pi i t).

$$\frac{8\,\Gamma\!\left(\dfrac{3}{4}\right)^{12}}{\pi^{3}}$$

2.95832521157700345043607659316706

1,24,312,2912,21816,139152,783328,3986112,18650424,81251896,332798544,
1291339296,4776117216,16922753616,57683178432,189821722688,604884735288,
1871370360240,5633654421720,16535803556064,47405095227984,
132942579098368,365211946954656

---

A319822 Number of solutions to x^2 + 2*y^2 + 5*z^2 + 5*w^2 = n.

$$\frac{\Gamma\left(\frac{5}{8}\right)^4 (3+2\sqrt{2})(5+\sqrt{5})^3 \sqrt{2+\sqrt{2}} \sqrt{5}}{8000\, \pi\, \Gamma\left(\frac{7}{8}\right)^4}$$

1.0904931781796092481506623925319

1,2,2,4,2,4,12,8,18,14,4,28,12,24,32,0,34,20,14,28,4,32,44,40,28,10,40,56,64,72,8,48,
66,24,68,8,46,88,60,32,4,52,64,116,76,12,64,72,60,82,26,72,104,104,88,8,112,56,
136,140,8,136,96,72,98,16,72,132

---

A320049 Expansion of (psi(x) / phi(x))^6 in powers of x where phi(), psi() are Ramanujan theta functions.

$$\frac{e^{\frac{3\pi}{4}} 2^{1/4}}{16}$$

0.78418725859165474593089953428806

1,-6,27,-98,309,-882,2330,-5784,13644,-30826,67107,-141444,289746,-578646,1129527,
-2159774,4052721,-7474806,13569463,-24274716,42838245,-74644794,128533884,
-218881098,368859591,-615513678,1017596115,-1667593666,2710062756,
-4369417452

---

A320050 Expansion of (psi(x) / phi(x))^7 in powers of x where phi(), psi() are Ramanujan theta functions.

$$\frac{e^{\frac{7\pi}{8}} 2^{5/8}}{32}$$

0.75304872679352721329186442864644

1,-7,35,-140,483,-1498,4277,-11425,28889,-69734,161735,-362271,786877,-1662927,
3428770,-6913760,13660346,-26492361,50504755,-94766875,175221109,
-319564227,575387295,-1023624280,1800577849,-3133695747,5399228149,
-9214458260,15584195428

---

A320069 Expansion of 1/(theta_3(q) * theta_3(q^2)), where theta_3() is the Jacobi theta function.

$$\frac{4\sqrt{\pi}\, \Gamma\left(\frac{7}{8}\right)^2 \sqrt{2-\sqrt{2}}}{\Gamma\left(\frac{5}{8}\right)^2 (2+\sqrt{2})}$$

0.91701683496666156658773258797020

1,-2,2,-4,10,-16,20,-32,58,-86,112,-164,260,-368,480,-672,986,-1348,1750,-2372,3312,
-4416,5684,-7520,10148,-13266,16912,-21960,28896,-37168,46944,-60032,77466,
-98312,123076,-155392,197422,-247696,307540,-384096,481776,-598500

---

A320070 Expansion of 1/(theta_3(q) * theta_3(q^2) * theta_3(q^3)), where theta_3() is the Jacobi theta function.

$$\frac{32 \pi^{15/4} \sqrt{3}\, \Gamma\!\left(\frac{7}{8}\right)^{8} \sqrt{2-\sqrt{2}}}{9\, \Gamma\!\left(\frac{5}{8}\right)^{8} \Gamma\!\left(\frac{7}{12}\right)^{3} \Gamma\!\left(\frac{2}{3}\right) \Gamma\!\left(\frac{11}{12}\right)^{2} (1+\sqrt{3})\left(\frac{17}{12}+\sqrt{2}\right)}$$

0.916868853218401561965239617 18355

1,-2,2,-6,14,-20,32,-60,98,-150,232,-360,558,-828,1196,-1776,2614,-3700,5238,-7480,
10516,-14592,20180,-27832,38216,-51970,70184,-94842,127612,-170140,226164,
-300324,396754,-521520,683484,-893432,1164330,-1511188,1954756,-2524188

---

A320124 Number of integer solutions to a^2 + b^2 + 2*c^2 + 3*d^2 = n.

$$\frac{2^{3/4}\, \Gamma\!\left(\frac{2}{3}\right) \Gamma\!\left(\frac{5}{8}\right)^{5} \Gamma\!\left(\frac{7}{12}\right) (2+\sqrt{2})^{3} (1+\sqrt{3})}{128\, \pi^{2}\, \Gamma\!\left(\frac{7}{8}\right)^{5}}$$

1.18494025334124335490988 89834320

1,4,6,10,20,20,24,40,22,28,56,20,50,80,28,80,84,32,78,80,68,100,120,80,88,124,56,82,
136,100,140,200,86,80,192,72,140,240,120,200,248,80,112,176,100,260,224,160,210,
172,186,128,272,180,240,400,124,200,280,116

---

A320126 Number of integer solutions to a^2 + b^2 + 2*c^2 + 5*d^2 = n.

$$\frac{\sqrt{2}\, 5^{3/4}\, \Gamma\!\left(\frac{5}{8}\right)^{4} (3+2\sqrt{2})(5-\sqrt{5})^{3/2} (\sqrt{5}+1)^{3} \sqrt{2+\sqrt{2}}}{6400\, \pi\, \Gamma\!\left(\frac{7}{8}\right)^{4}}$$

1.18474939307740355720555 91445492

1,4,6,8,12,10,16,28,22,36,40,24,56,36,24,64,28,64,78,24,58,48,68,92,80,92,72,112,92,
48,104,96,118,176,64,48,124,84,148,160,104,120,176,120,72,146,88,204,216,124,
126,96,180,148,224,188,120,304,216,120,224,96

---

A320138 Number of integer solutions to a^2 + 2*b^2 + 3*c^2 + 3*d^2 = n.

$$\frac{\Gamma\left(\frac{2}{3}\right)^2 \Gamma\left(\frac{5}{8}\right)^6 \Gamma\left(\frac{7}{12}\right)^2 (7\sqrt{2}+10)(2+\sqrt{3})\sqrt{2}\sqrt{2+\sqrt{2}}}{128\pi^3 \Gamma\left(\frac{7}{8}\right)^6}$$

1.0908445581114080430744973274774

1,2,2,8,10,8,24,16,10,38,8,12,48,8,32,64,26,36,70,28,24,80,28,48,96,42,40,76,48,24,112, 64,58,160,68,32,126,56,44,192,56,84,176,44,60,88,80,96,208,114,74,176,72,72,172, 80,112,288,88,76,224,72,112,304,90,96

---

A320140 Number of integer solutions to a^2 + 2*b^2 + 3*c^2 + 5*d^2 = n.

$$-\frac{128(-3+\sqrt{3})\,10^{1/4}(2-\sqrt{2})^{9/2}\left(\frac{5}{2}+\sqrt{5}\right)\pi^4 \Gamma\left(\frac{7}{8}\right)^9 (5-\sqrt{5})^{3/2}}{15\sqrt{2+\sqrt{2}}(\sqrt{5}+1)^3 \Gamma\left(\frac{11}{12}\right)^2 \Gamma\left(\frac{7}{12}\right)^3 \Gamma\left(\frac{5}{8}\right)^9 \Gamma\left(\frac{2}{3}\right)}$$

0.91686857687104464057169857981772

1,2,2,6,6,6,16,8,14,26,8,32,26,8,40,16,22,40,22,32,46,40,24,48,40,42,72,50,32,64,56,28, 74,48,60,112,78,24,72,76,40,144,48,48,120,50,52,48,70,98,150,128,40,84,128,52,176, 120,56,208,96,72,92,72,102,192,156

---

A320147 Number of integer solutions to a^2 + b^2 + 3*c^2 + 5*d^2 = n.

$$\frac{\Gamma\left(\frac{7}{10}\right)\Gamma\left(\frac{7}{12}\right)\Gamma\left(\frac{2}{3}\right)2^{1/10}5^{1/4}\Gamma\left(\frac{3}{5}\right)^2(\sqrt{5}+1)^2(1+\sqrt{3})(5+\sqrt{5})}{320\sqrt{\pi}\,\Gamma\left(\frac{3}{4}\right)^5 \Gamma\left(\frac{9}{10}\right)}$$

1.180531460666086416195687473151

1,4,4,2,12,18,8,16,24,28,40,8,26,72,16,16,44,44,68,24,34,80,72,28,40,124,40,50,112,56, 80,40,76,144,120,32,84,216,24,40,136,80,160,88,56,154,88,28,158,228,100,48,216, 172,80,104,80,300,280,40,112,248,120,112

---

A320149 Number of integer solutions to a^2 + 2*b^2 + 2*c^2 + 2*d^2 = n.

$$-\frac{\Gamma\left(\frac{5}{8}\right)^4 (3+2\sqrt{2})\sqrt{2+\sqrt{2}}}{16\pi\,\Gamma\left(\frac{7}{8}\right)^4 (\sqrt{2}-2)}$$

1.0986534618873295727942351752354

1,2,6,12,14,24,20,16,30,14,40,60,36,72,48,16,62,36,42,108,72,96,100,48,68,42,120,120, 112,168,48,64,126,40,108,192,98,216,180,48,136,84,160,252,180,168,144,96,132,

114,126,216,216,312,200,80,240,72,280,348,112

---

A320150 Number of integer solutions to a^2 + 2*b^2 + 2*c^2 + 3*d^2 = n.

$$\frac{3^{3/4}\,\Gamma\!\left(\frac{2}{3}\right) 2^{1/4}\,\Gamma\!\left(\frac{5}{8}\right)^3 (2+\sqrt{2})^{3/2}}{48\,\pi\,\Gamma\!\left(\frac{7}{8}\right)^3 \Gamma\!\left(\frac{11}{12}\right) (\sqrt{2}-1)(\sqrt{3}-1)}$$

1.0947420473107136053688100330414

1,2,4,10,10,16,24,12,28,26,8,48,30,28,72,32,34,64,28,36,80,60,72,96,72,42,56,82,36,
112,120,60,124,96,32,96,130,76,216,140,56,160,48,84,144,112,144,192,150,86,84,
128,140,208,240,96,216,180,56,240,96,124,360

---

A320151 Number of integer solutions to a^2 + 2*b^2 + 2*c^2 + 5*d^2 = n.

$$\frac{5^{3/4}\,\Gamma\!\left(\frac{5}{8}\right)^4 (3+2\sqrt{2})\sqrt{2}\,(5-\sqrt{5})^{3/2}(\sqrt{5}+1)^3(4+2\sqrt{2})}{25600\,\pi\,\Gamma\!\left(\frac{7}{8}\right)^4}$$

1.0945657154193823122175866877714

1,2,4,8,6,10,12,8,28,22,24,40,8,32,36,16,54,32,28,56,26,32,48,40,84,74,72,48,24,92,52,
96,92,32,96,80,42,64,80,80,168,124,48,72,72,94,132,104,72,126,124,96,48,96,120,
168,252,96,120,168,48,196,128,88,246

---

A320239 Expansion of theta_3(q) * theta_3(q^3) * theta_3(q^5), where theta_3() is the Jacobi theta function.

$$-\frac{\Gamma\!\left(\frac{11}{12}\right)^2 \Gamma\!\left(\frac{7}{12}\right)^3 \sqrt{3}\,\Gamma\!\left(\frac{2}{3}\right) \pi^{1/4} 5^{3/4} (1+\sqrt{3})^3 (-2+\sqrt{3})(5-\sqrt{5})^{3/2}(\sqrt{5}+1)^3}{6400\,\Gamma\!\left(\frac{3}{4}\right)^8}$$

1.0866104882516378755410636489104

1,2,0,2,6,2,4,4,4,14,0,0,14,4,4,0,6,12,8,4,2,20,0,4,20,2,8,10,12,4,4,4,16,32,0,0,26,4,0,
12,0,20,8,4,8,6,4,4,42,18,0,8,20,12,16,0,12,48,8,8,0,16,8,12,14,0,16,4,20,24,4,0,36,
28,0,2,20,8,8,4,6

---

A330373 Sum of all parts of all self-conjugate partitions of n.

$$\frac{e^{-\frac{\pi}{24}}\,2^{1/4}}{24}$$

0.043470760943528625521994820367317

0,1,0,3,4,5,6,7,16,18,20,22,36,39,42,60,80,85,90,114,140,168,176,207,264,300,312,378,
448,493,540,620,736,825,884,1015,1188,1295,1406,1599,1840,2009,2184,2451,2772,
3060,3312,3666,4176,4557,4900,5457,6084,6625,7182,7920,8792,9576,10324,11328,
12540

---

### A347801 Expansion of ( Sum_{k>=0} k^2 * q^(k^2) )^2.

$$\frac{1}{64\,\pi^{3/2}\,\Gamma\!\left(\frac{3}{4}\right)^2}$$

0.0018686485405221383684606569723056

0,0,1,0,0,8,0,0,16,0,18,0,0,72,0,0,0,32,81,0,128,0,0,0,0,288,50,0,0,200,0,0,256,0,450,0,
0,72,0,0,288,800,0,0,0,648,0,0,0,0,723,0,1152,392,0,0,0,882,0,0,1800,0,0,0,1696,0,
0,512,0,0,0,1296,1152,2450,0,0,0,0,0,2048,0,162,0,0,4176,0,0,0,3200,1458

---

### A347802 Expansion of ( Sum_{k>=0} k^2 * q^(k^2) )^3.

$$\frac{1}{512\,\pi^{9/4}\,\Gamma\!\left(\frac{3}{4}\right)^3}$$

0.000080777691771496507346868329874754

0,0,0,1,0,0,12,0,0,48,0,27,64,0,216,0,0,432,48,243,0,384,972,0,768,0,864,804,0,3456,
600,0,0,1968,3888,1350,3072,0,5508,0,0,7776,2400,6075,1728,9600,1944,0,4096,
7776,21600,2022,0,3456,17424,0,13824,21552,0,19521,0,31104,15984,0,0,21600,
34896,11907

---

### A347803 Expansion of ( Sum_{k>=0} k^2 * q^(k^2) )^4.

$$\frac{1}{4096\,\pi^3\,\Gamma\!\left(\frac{3}{4}\right)^4}$$

$3.4918473679955178000766058732981 \times 10^{-6}$

0,0,0,0,1,0,0,16,0,0,96,0,36,256,0,432,256,0,1728,64,486,2304,768,3888,0,3072,7776,
1728,7112,0,13824,12864,0,27648,6336,15552,9261,18688,62208,21744,24576,0,
72576,51456,24300,117504,38400,101088,9216,93184,155520,86400,142382,62208,
352512,67344,0,202752,286176

---

### A350642 Expansion of Product_{k>=1} (1-q^(2*k))/(1-q^k)^4.

$$\frac{2 e^{-\frac{\pi}{12}} \Gamma\left(\frac{3}{4}\right)^3}{\pi^{3/4}}$$

1.2003860731402909414643714570620

1,4,13,36,90,208,455,948,1901,3688,6955,12792,23019,40612,70395,120072,201822, 334684,548158,887500,1421602,2254460,3541928,5515900,8519173,13055208, 19859113,29998024,45012751,67116436,99472320,146580028,214811311,313149460

---

A350643 Expansion of Product_{k>=1} (1-q^(2*k))^2/(1-q^k)^7.

$$\frac{2 e^{-\frac{\pi}{8}} 2^{5/8} \Gamma\left(\frac{3}{4}\right)^5}{\pi^{5/4}}$$

1.3759680430559075787477955224413

1,7,33,126,419,1260,3509,9185,22842,54395,124784,277059,597644,1256341,2580363, 5189185,10236710,19840410,37832553,71060190,131610897,240585292,434431132, 775483785,1369359198,2393425484,4143057525,7106240582,12083072562, 20375932566

---

A350644 Expansion of Product_{k>=1} (1-q^(2*k))^3/(1-q^k)^10.

$$\frac{4 e^{-\frac{\pi}{6}} 2^{1/4} \Gamma\left(\frac{3}{4}\right)^7}{\pi^{7/4}}$$

1.5772325181345541914860946420383

1,10,62,300,1235,4522,15130,47084,137990,384370,1024760,2629380,6521693, 15693180,36745810,83935920,187441365,409981826,879717860,1854439520, 3845126929,7850815860,15799770260,31368976420,61490409175,119092108534, 228039325630

---

A360191 G.f. 1 / Product_{n>=1} (1 - x^n)^3 * (1 - x^(2*n-1))^2.

$$\frac{e^{-\frac{\pi}{24}} 2^{7/8} \Gamma\left(\frac{3}{4}\right)^3}{\pi^{3/4}}$$

1.2547038002583000343985678035316

1,5,18,55,149,371,867,1923,4086,8374,16634,32152,60669,112041,202943,361200, 632647,1091917,1859225,3126242,5195715,8541624,13899866,22404091,35787815, 56683294,89061028,138872410,214984454,330532633,504869316,766357010,

*1156355165*

---

*A361535 Expansion of g.f. 1 / Product_{n>=1} ((1 - x^n)^6 * (1 - x^(2*n-1))^4).*

$$\frac{2\,e^{-\frac{\pi}{12}}\,2^{3/4}\,\Gamma\!\left(\frac{3}{4}\right)^6}{\pi^{3/2}}$$

1.57428162638262006946674638306677

*1,10,61,290,1172,4212,13833,42262,121625,332764,871641,2197936,5359005, 12679730,29200593,65617892,144189054,310400110,655669910,1360910666, 2779007594,5589070978,11081585154,21679798590,41883282555,79958881544, 150943109191,281926365224*

---

*A385520 Expansion of Product_{k>0} ((1 - q^(2*k))*(1 - q^(6*k))^3)/((1 - q^k)*(1 - q^(3*k))*(1 - q^(4*k))*(1 - q^(12*k))).*

$$\frac{3^{5/12}\,\Gamma\!\left(\frac{2}{3}\right)^{1/3}\,\Gamma\!\left(\frac{3}{4}\right)^{7/3}\,(\sqrt{2}\,(1+\sqrt{3}))^{4/3}\,(\sqrt{3}-1)}{6\,\pi^{1/6}\,\Gamma\!\left(\frac{11}{12}\right)^{4/3}\,\Gamma\!\left(\frac{7}{12}\right)}$$

1.04533820420055945597529278219l8

*1,1,1,3,4,5,6,9,13,16,20,27,36,44,54,69,88,107,130,162,200,240,288,351,426,507,602, 723,864,1019,1200,1422,1681,1968,2300,2700,3160,3674,4266,4965,5768,6665,7692, 8892,10260,11792,13536,15552,17844,20407*

---

(1)